**Title:** Control System Design Using Finite Laplace Transform Theory

**Author:** Subhendu Das, CCSI, California, USA

Author email: subhendu11das@gmail.com

*Abstract* – The Laplace transform theory violates a very fundamental requirement of all engineering systems. We show that this theory assumes that all signals must exist over infinite time interval. Since in engineering this infinite time assumption is not meaningful and feasible, this paper presents a design for linear control systems using the well known theory of Finite Laplace transform (FLT). The major contributions of this paper can be listed as: (a) A design principle for linear control systems using FLT, (b) A numerical inversion method for the FLT with examples, (c) A proof that the FLT does not satisfy the convolution theorem as normally required in engineering design and analysis, and (d) An observation that the FLT is conceptually similar to the analog equivalent of the Finite Impulse Response (FIR) digital filter.

I. Introduction

Most of our engineering systems run over finite time. Consider the example of a robot arm, picking up an item from one place and dropping it in another place and repeating the process in, say, less than a second of time. Similarly a digital communication receiver system, receives an electrical signal of microsecond duration, for example, representing the data, extracts the data from the signal, sends it to the output, and then goes back to repeat the process. This is the general nature of our technology today. Although many of our engineering systems run over infinite time, like global positioning system receivers in our cell phones, traffic light signaling systems at the street intersections etc., but if you examine the internal implementation details you will find that the repeating processes of finite time are everywhere.

Our software is running under operating systems which are also nothing but finite state machines. In addition, the software system is very dynamic also. It allows us to periodically monitor the system, stop the processes, examine its status, and restart in a different way if required. Thus the modern engineering is continuously interruptive not only from external stimuli, but also from internal status. Thus the finite time is an integral part of today's engineering systems. The continuous analog systems do not exist anymore. In all the examples mentioned above, the systems take different finite time intervals for different activities, but in this paper for simplicity and without loss of generality, we will assume that all finite time intervals are same and fixed.

The main objective of this research is to emphasize that engineering design should be treated as a finite time design problem. Applying an infinite time theory, like Laplace Transform, to a finite time problem appears to be illogical. A finite time approach will help to keep it consistent with engineering requirements and will help to improve the product quality. The ideas presented in this paper are very simple but we still discuss it in details to make the paper self contained, coherent, and readable for most engineers in the industry. It should be recognized that the Infinite Laplace Transform theory (ILT) is still a predominant tool for engineers in both commercial and military industries. This paper is not for presenting or comparing other better methods available for the design of control systems. It's only focus is to bring out the fact that ILT is not the correct tool for engineering designs and also to show that, if we still need to stay in the Laplace domain, then we should use FLT. As we will see, that the FLT theory is more realistic, simpler to understand, and easier to apply.

The paper is organized along the following lines. In the first three sections we introduce the FLT, some of its basic properties that are relevant to control systems, and discuss inverse FLT. In section four we present the numerical inversion method for the FLT. Section five presents an application to simultaneous differential equations showing how to resolve the boundary conditions when we use FLT. In section six we show that the FLT does not satisfy the standard convolution theorem. Then in section seven we present the FLT based design method for the control system. Finally we briefly talk about the finite time issues related to the Fourier transform theory in section eight. During this presentation we point out in many places and in many different ways what will happen if we use ILT for finite time systems. We also clearly highlight the contributions made by this paper in the FLT theory as well as in the control engineering at all relevant locations throughout the paper.

## II. Finite Laplace Transform

The standard Laplace transform is defined by the following expression [1]:

$$F(s) = \int_0^\infty e^{-st} f(t) dt \tag{1}$$

Since the upper limit is infinity, this definition will be referred to as the Infinite Laplace Transform (ILT). Because of this infinite limit, the integral (1) requires a convergence or boundedness condition. To ensure convergence it is assumed that the growth of the function f(t) is restricted by the following constraint:

$$|f(t)| \leq Me^{\alpha t}, \quad 0 < M < \infty \tag{2}$$

We will assume that all our functions are continuous since we are considering only engineering problems. If a continuous function satisfies (2) then it is called a function of exponential order α. It can be shown [1, p. 13] that the integral (1) converges for all $\mathcal{R}(s) > \alpha$. Here $\mathcal{R}(s)$ means the real part of the complex variable s. This region is called the region of convergence. We now present the well known [2] basic theories of FLT with some examples.

Along the definition of ILT (1) the Finite Laplace Transform (FLT) has been defined as [2, pp 283-294]

$$\mathcal{L}_T[f] = F_T(s) = \int_0^T e^{-st} f(t) dt \; , \; 0 < T < \infty \tag{3}$$

From now on we will use F(s) to denote FLT of a continuous function f(t). In (3) t will be referred to as time, [0,T] as finite duration interval, s as a complex variable represented by x+iy, where x and y are real variables. We will assume that T is constant and represent a fixed time interval on which our engineering processes will repeat. Since the upper limit is finite, the integral (3) always exist. Thus the region of convergence of FLT is the entire complex plane. This is a significant difference of FLT compared to ILT of (1). Although this property is well known in the FLT theory, but it has not been emphasized and used in engineering. This paper will utilize the full potential and consequences of this analyticity property.

It should be clear that we are considering only continuous functions. In engineering applications there are no piecewise continuous functions. If you put an oscilloscope probe at any pin of a digital microprocessor then you will always find a continuous signal. Thus even a digital system does not use piecewise continuous functions. We will also not consider Stieltjes integrals, functions of bounded variations, etc., in this paper, for simplicity and because they are not relevant in engineering.

Let us consider two examples of FLT [2]. The finite duration step function f(t) is defined by

$$f(t) = \begin{cases} 1 & 0 \leq t \leq T \\ 0 & otherwise \end{cases} \tag{4}$$

Using the definition (3) we get

$$\mathcal{L}_T(1) = \int_0^T e^{-st} \cdot 1 \cdot dt$$

$$= \frac{1}{s} - \frac{1}{s}e^{-sT} \tag{5}$$

$$= \frac{1 - e^{-sT}}{s} \tag{6}$$

We can see from (5) that the FLT has the standard ILT term $\frac{1}{s}$ and it also has another expression involving $e^{-sT}$. We will see that this exponential term will always be there in all FLT expressions. They can be verified from (3) easily. If we use only the first part of (5), that is, the ILT part, then we will implicitly assume infinite time situation for our finite time problems and the Laplace model will not be correct for the engineering applications. Thus whenever we are using ILT only expressions, we are automatically assuming an infinite time model. It is also worth noting at this time that (6) does not have a pole at s=0, because the numerator also goes to zero at s=0. This property of FLT, 0/0 form, is a very distinguishing feature. This property also has not been clearly stated in the literature. We will use it extensively in all the theories presented here. This is of course a consequence of analytic property of the FLT theory. We will see that this will hold in all cases, and will give a theoretical proof later in this paper to show that the FLT is indeed an analytic function over the entire complex plane. Thus when we apply FLT we do not get any poles in the transfer function. Thus treating finite time problems with poles, as in the ILT case, is very artificial, the real system actually does not have poles. In this paper we will use analytical models like (6) instead of the first part of (5).

We will not derive many examples here because of lack of space and because they are all available in the literature [2]. Let us examine the FLT of an exponential function.

$$\mathcal{L}_T(e^{at}) = \int_0^T e^{at} e^{-st} dt = \int_0^T e^{-(s-a)t} dt$$

$$= \left.\frac{e^{-(s-a)t}}{-(s-a)}\right|_0^T = \frac{1}{s-a} - \frac{1}{s-a} e^{-(s-a)T} \tag{7}$$

$$= \frac{1 - e^{-(s-a)T}}{s-a} \tag{8}$$

Again we see that (7) has the ILT expression in addition to the exponential term. Also at s=a the expression (8) takes the 0/0 form showing no poles at s=a, unlike ILT. For later references we also record the FLT pairs for sine and cosine functions.

$$\sin wt \leftrightarrow \frac{w}{s^2+w^2} - e^{-sT}\frac{w}{s^2+w^2}\cos wT - e^{-sT}\frac{s}{s^2+w^2}\sin wT \tag{9}$$

$$\cos wt \leftrightarrow \frac{s}{s^2+w^2} - e^{-sT}\frac{s}{s^2+w^2}\cos wT + e^{-sT}\frac{w}{s^2+w^2}\sin wT \tag{10}$$

The above examples show that the FLT expressions are more complex to deal with, and may be because of this reason, the infinite limit was considered originally in the definition of (1). Note that if we let T go to infinity all the FLT expressions become standard ILT expressions.

From the fundamental theorem of calculus [3, p. 142] we know that if f(t) is continuous on [a,b] and F(t) is the antiderivative of f(t) defined by

$$F(t) = \int f(t)\,dt \quad \text{then} \quad \int_a^b f(t)dt = F(b) - F(a)$$

Using the above fundamental theorem if we define

$$F(s,t) = \int e^{-st} f(t)dt \quad \text{then it is clear that}$$

FLT:  $F(s,T) - F(s,0)$ \hfill (11)

ILT:  $-F(s,0)$ \hfill (12)

It is easily seen from (1) and (2) that the ILT has been defined in such a way that $F(s,\infty)$ is zero and we only have $-F(s,0)$ in (12), the same part that we always get in the FLT (11).

It seems that the FLT theory was first introduced by Dunn [4]. Later another paper was published by Debnath [5], which is also available as a chapter in the book by Debnath [2]. The objectives of both papers were mathematical in nature, to extend the ILT theory to cover a larger class of Laplace transformable functions. For example using FLT we can take transform of $\exp(t^2)$ which is not a function of exponential order. Those two papers and the book do not present any engineering applications. The motivation of this paper though, is driven by the finite time requirement of all engineering applications, which somehow appears to have never been pointed out in the Laplace literature. A detailed examination of the embedded hardware, the software, and the systems will however, reveal this finite time feature [6, pp. 73-88]. The examples given at the beginning of the introduction section are very convincing about the finite time nature of modern engineering.

## III. FLT Properties

In this section we discuss some of the important properties of the Finite Laplace Transform (FLT) theory relevant to control system engineering. Some of the results, comments, and the explanations provided are new and did not appear in the literature before.

*3.1 Convergence*

We have mentioned that the FLT definition (3) always converges. This is mainly because FLT integration is taken over finite duration interval and we have considered continuous functions only. Thus the region of convergence of the FLT is the entire complex plane. This essentially means that all FLT expressions are valid for all s in the complex plain. And the complex function F(s) never goes to infinity, that is, there are no poles or singularities in the FLT expressions. Thus if we use ILT expressions for finite time problems we cannot get the correct features of the model as mentioned before.

*3.2 Analyticity*

It was briefly mentioned that F(s) is an analytic function. The proof is quite simple but it has not appeared in [2] or [4]. Yamamoto [7] points to the Paley-Wiener theorem for the proof. A function F(s) of complex variable s is analytic at a point $s_0$ if F(s) has a derivative at some neighborhood of $s_0$. It should be pointed out that $s_0$ is included in that neighborhood. A function F(s) is called an entire function if it has derivatives at each nonzero point in the finite complex plane [3, p. 73]. We show that the FLT is indeed an entire function. Another mathematical term for analytic function is holomorphic. We state the following new theorem and its proof.

***Theorem***: The complex FLT function F(s) satisfies the following conditions:

(a) $F(s) = u(x,y) + iv(x,y)$ is defined over the entire complex plane,

(b) The first order partial derivatives of u and v exist at all points s=x+iy in the plane, and

(c) The partial derivatives are continuous and satisfy the Cauchy-Riemann equations

$$u_x = v_y, \quad u_y = -v_x,$$ at all points in the plane.

Therefore $F'(s)$ exists at all points in the complex plane and F(s) is an entire function.

The proof is quite straight forward and follows the lines of ILT similar to [1, pp. 124-125]. For any s=x+iy in the finite plane we can write from (3)

$$F(s) = \int_0^T e^{-st} f(t)dt = \int_0^T e^{-(x+iy)t} f(t)dt$$

$$= \int_0^T e^{-xt} \cos yt \, f(t)dt + i \int_0^T e^{-xt} \sin yt \, f(t)dt$$

$$= u(x, y) + iv(x, y) \tag{13}$$

It is clear from (13) that u and v are defined and exist since the functions are continuous and the integration is over finite duration. Also,

$$u_x = \frac{\partial}{\partial x}\left[\int_0^T e^{-xt} \cos yt \, f(t)dt\right] = \int_0^T \frac{\partial}{\partial x}[e^{-xt} \cos yt \, f(t)dt] \tag{14}$$

We can see that the partial derivative in (14) can be moved from outside the integral to inside without any problem because of finite duration and the continuity of all functions involved. Thus there are no absolute or uniform convergence issues to be considered like in ILT case. Hence we can write from the last expression

$$u_x = \int_0^T \frac{\partial}{\partial x}[e^{-xt} \cos yt \, f(t)dt] = \int_0^T -te^{-xt} \cos yt \, f(t)dt \tag{15}$$

Similarly we can show that

$$v_y = \int_0^T \frac{\partial}{\partial y}[-e^{-xt} \sin yt \, f(t)dt] = \int_0^T -te^{-xt} \cos yt \, f(t)dt \tag{16}$$

The above two expressions, (15) and (16), show that $u_x = v_y$.

In the similar way we can also show that $u_y = -v_x$.

The two equality conditions on partial derivatives, called Cauchy-Riemann conditions [3, p. 66], indicates that the derivative of F(s) exists for all s in the entire plane and thus by definition F(s) is an entire function. This completes the proof.

The above proof is quite simple. The proof of Paley Weiner theorem can be found in Rudin [8, pp. 180-185]. The main idea behind this analyticity result is that if we extend the ILT to finite time problems then there will be no poles in the complex plane. As a consequence of this, the residue method cannot be used to find the inverse of the FLT functions. For the contour integral will be zero around any closed contour. However, many authors, [4], [5], and others, have treated the two parts of the FLT expression (11) independently and used contour integration to find inverse. In this paper we consider the composite expression (11) together as an analytic function. Somehow the literature has ignored this composite approach and used the classical residue based methods. This paper treats the composite function as a whole and thus exploits the full potential of the analyticity.

*3.3 Taylor Series*

One of the fascinating properties of functions of complex variables is that if it is analytic then it has all the derivatives. Since the FLT F(s) is analytic in the entire complex plane it has the following Taylor series expansion for F(s) around the origin and is also valid over the entire complex plane [3, pp. 189-192].

$$F(s) = \sum_{n=0}^{\infty} \frac{F^{(n)}(0)}{n!} s^n \tag{17}$$

Therefore F(s) can be approximately expressed by the finite series:

$$F(s) \approx a_0 + a_1 s^1 + a_2 s^2 + \cdots + a_N s^N \tag{18}$$

Interestingly the expression (17) is not valid for ILT since ILT has poles in the left half plane. When there are singularities or poles, F(s) can be expressed as infinite series, called Laurent series, which has both positive and negative powers of the complex variable s. We will show later how expression (18) can be used to numerically invert the FLT of any continuous function over the finite duration. This Taylor series approach is also new in this paper. The FLT has not been treated before in the literature using the Taylor series, particularly the way we have done it here, as far as we know.

*3.4 FIR Filter*

A Finite Impulse Response (FIR) filter is a very popular digital signal processing filter. It uses the sampled values of the continuous time input signal and is defined as a linear combination of present and previous values of the input signal. It does not use the previous outputs from the filter. An Infinite Impulse Response (IIR) filter on the other hand uses both previous input and output values of the signal. A FIR filter can be expressed as [9, pp. 148-158]

$$H(z) = \sum_{k=0}^{N-1} h(k)\, z^{-k} \qquad (19)$$

Here H(z) is the transfer function of the FIR filter, z is the variable for the Z-Transform, and h(k) is the k-th time sample of the impulse response of the continuous time signal.

Thus H(z) is a polynomial in z and has no poles in the Z-plane, except possibly at z=0, which is equivalent to the negative infinity in the Laplace plane. Because of this reason an FIR filter is always stable. The FLT shares these two properties with the FIR filter. It is worth noting also that the expressions (18) and (19) are very similar. Thus we can say that FLT is the continuous time analog equivalent of the discrete time FIR digital filter.

*3.5 Stability*

We make some general comments here on the classical theory of the ILT and in the light of the finite time concepts of engineering. The stability concept is related to the behavior of a function f(t) as t approaches to infinity. Since in FLT and in all engineering problems infinite time is not used and is not meaningful, the usual stability concept is also not meaningful. The complex FLT function F(s) is analytic and therefore does not go to infinity. The FLT does not have any singularity at any point in the entire plane, the corresponding time function therefore also cannot go to infinity. That is because all singularities of the ILT are related to the poles and they translate to exponential functions.

Because of the same reason the final value theorem [1, p. 89]

$$\lim_{t \to \infty} f(t) = \lim_{s \to 0} F(s)$$

is also not applicable in engineering. In the above expression F(s) is the ILT. Clearly t going to infinity cannot be meaningful in any engineering applications because we have mentioned in many different ways that all engineering applications are performed over finite time only. Thus stability theory helps to understand the behavior of systems in an analytical frame work.

Since the functions do not have time to stabilize in engineering, the notion of frequency response [10] will also not be meaningful for FLT theory. The working system will not see the stable steady state response during time T for a sinusoidal input frequency signal. Therefore the gain, which is the ratio of the magnitude of the steady values of the output to that of the input will not be defined. The same is true for the phase response, which compares the phase angle differences between the input, and the output sinusoidal response at steady state. These concepts are related to ILT and associated poles, and are not meaningful for FLT and in engineering. In FLT engineering the time period T is very small and the transients may not settle and thus the concept of steady state cannot be feasible or meaningful during that time period.

An interesting related concept is that if we go deeper to the details, and look at the register of a microprocessor, of an engineering real time embedded system, that executes at nanoseconds clock cycles at a periodic rate of few microseconds, you will not see any steady values there. There are many software engineering development tools that allow us to examine these registers at that time resolution. Thus the concept of stability at the detailed level is non-existent.

### 3.6   *Inverse FLT*

The FLT theory will be incomplete if we cannot find the inverse. For the Infinite Laplace Transform (ILT) theory the inverse is defined by the following relation [1, pp. 152-157]:

$$f(t) = \frac{1}{2\pi i} \int_{x-i\infty}^{x+i\infty} e^{ts} F(s) ds \qquad (20)$$

Here x > α, and α is defined in (2). It can be shown that (20) is equivalent to the following contour integral:

$$f(t) = \frac{1}{2\pi i} \oint_C e^{ts} F(s) ds \tag{21}$$

Here the contour C is taken appropriately to cover all the singularities on the left of the vertical line, called Bromwich line, passing through the point α. Normally (21) is computed using the relation

$$\frac{1}{2\pi i} \oint_C e^{ts} F(s) ds = \sum_{k=1}^{N} \lim_{s \to z_k} (s - z_k) e^{ts} F(s) \tag{22}$$

Here $z_k$ is a pole of F(s). The approach based on (22) is known as the residue method. Since FLT does not have any poles, the method in (22) cannot be used to find the inverse of the FLT functions. In fact we can see what happens if we try to use (22) for the inverse FLT of the exponential function:

$$\sum_{k=1}^{1} \lim_{s \to a} (s - a) e^{ts} F(s) = \lim_{s \to a} (s - a) e^{ts} \left[ \frac{1 - e^{-(s-a)T}}{s - a} \right]$$

$$= e^{at}(1 - e^{-0T}) = 0$$

Thus one of the very useful methods for ILT is no longer valid for FLT. However, as mentioned before, many authors [4,5] have treated the two terms in expression (11) separately using separate contours. In the following section we present a numerical method for the inversion of the FLT that uses the analyticity property of (11).

## VI    Numerical Inversion

There are many numerical methods for inverting the ILT functions. Probably Bellman was the first person to introduce the numerical approach concept in Laplace transform theory. Bellman [11, pp. 135-155] uses positive real integers for s before integrating the ILT equation

(1) as real and imaginary integrals. [12,13] gives a summary of recent most popular numerical methods, including Post-Widder and Talbot's methods. Some of them use deformed Bromwich contours. Numerical inversion of ILT is also an ill-posed problem as pointed out in [14].

We present here a numerical inversion method that takes advantage of the analyticity property of the FLT functions. This method is new and does not depend on the previous methods, which are mainly based on ILT concepts. This method cannot be used for ILT, for it is based on the Taylor's series and the ILT does not have any Taylor series expansion. Thus this approach is unique to FLT; however it has some similarity with the FIR filter methods. The objective is to find the time function f(t) from the known FLT function F(s). Since F(s) is analytic, it has the Taylor series:

$$F(s) = \int_0^T e^{-st} f(t)dt = a_0 + a_1 s + a_2 s^2 + a_3 s^3 \cdots \tag{23}$$

Expanding the exponential function $e^{-st}$ under the integral sign, ($e^{-st}$ is an analytic function), and equating the like powers of s in (23) it can be shown that

$$a_n = (-1)^n \frac{1}{n!} \int_0^T t^n f(t)dt \tag{24}$$

Thus the problem reduces to finding f(t) in (24) given all the coefficients $\{a_n, n = 0,1,\cdots,N\}$. It is possible to solve (24) in many different ways. We will use the simplest method to demonstrate the feasibility of the inverse FLT concept. We expand f(t) in Taylor series, assuming that f(t) has all the derivatives, which may not be practical even in engineering.

$$f(t) = b_0 + b_1 t + b_2 t^2 + b_3 t^3 \cdots \tag{25}$$

We can now substitute (25) in (24) and express each coefficient $a_n$ in terms of $b_n$ for n=1...N, to get

$$a_n = \sum_{k=1}^{N}(-1)^n \frac{1}{n!} \frac{1}{(n+k+1)} T^{n+k+1} b_k \tag{26}$$

For N=4, as an example, using matrix notation we can express equation (26) in more readable form as in (27). We can now solve the matrix equation (27) by inverting the square matrix to find the b-coefficients $\{b_k\}$. Clearly, the matrix is non-singular, because the columns are independent. These b-coefficients can then be used in (25) to find the function f(t).

$$\begin{bmatrix} a_0 \\ a_1 \\ a_2 \\ a_3 \end{bmatrix} = \begin{bmatrix} T & \frac{T^2}{2} & \frac{T^3}{3} & \frac{T^4}{4} \\ -\frac{T^2}{2} & -\frac{T^3}{3} & -\frac{T^4}{4} & -\frac{T^5}{5} \\ \frac{1}{2!}\frac{T^3}{3} & \frac{1}{2!}\frac{T^4}{4} & \frac{1}{2!}\frac{T^5}{5} & \frac{1}{2!}\frac{T^6}{6} \\ -\frac{1}{3!}\frac{T^4}{4} & -\frac{1}{3!}\frac{T^5}{5} & -\frac{1}{3!}\frac{T^6}{6} & -\frac{1}{3!}\frac{T^7}{7} \end{bmatrix} \begin{bmatrix} b_0 \\ b_1 \\ b_2 \\ b_3 \end{bmatrix} \quad (27)$$

As an example consider the FLT expression (8) for the exponential function. The inverse was generated using Mathematica analysis tool. The Figure 1 shows how the inverse compares with the exact exponential function with four coefficients of the FLT Taylor series (23). For 10 FLT coefficients the true graph and the numerically reconstructed graph match exactly at the resolution of the paper. The numerical data are shown in the Table 1. The data were normalized, and therefore the last sample was always 1.0. Thus this Taylor series based approach appears to be very robust and satisfactory even for small number of terms. This does not normally happen even for FIR filters.

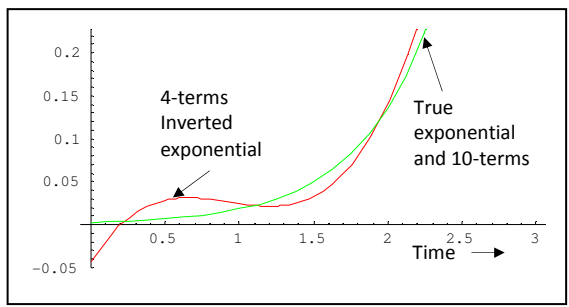

Fig.1: Inverse FLT using four coefficients

| True exponential data samples | | | | |
|---|---|---|---|---|
| 0.00451658 | 0.00822975 | 0.0149956 | 0.0273237 | 0.0497871 |
| 0.090718 | 0.165299 | 0.301194 | 0.548812 | 1.0 |
| Exponential from inverse FLT with four coefficients | | | | |
| 0.0133275 | 0.030779 | 0.0271752 | 0.0214159 | 0.0324014 |
| 0.0790315 | 0.180206 | 0.354826 | 0.621791 | 1.0 |
| Exponential from inverse FLT with ten coefficients | | | | |
| 0.00451518 | 0.00823095 | 0.0149953 | 0.0273231 | 0.0497887 |
| 0.090718 | 0.165299 | 0.301198 | 0.548813 | 1.0 |

Table 1: Numerical inversion of FLT

Taylor series approach has been used [15] for ILT inversion analytically. Their approach expands the function f(t) in Taylor series and evaluates over discrete values of s. In our paper we have expanded all functions, f(t), F(s), and $e^{-st}$ over all s. Also our approach is a numerical approach for FLT.

## V   Applications

Many standard results of ILT can be extended in a straight forward way using the methods shown in [1] for the case of FLT. These are for example, translation in t-domain, translation in s-domain, derivative of F(s) etc. We only present some results that will be required for our main objective and, although may be trivial, are not available in the literature. We also explain clearly, which is not apparent in the existing literature [2,4,5], the method of finding the boundary conditions using analyticity property of FLT. Consider the FLT of the first derivative with respect to time of a function

$$\mathcal{L}_T[f'(t)] = \int_0^T e^{-st} f'(t)\, dt$$

$$= e^{-st} f(t)\big|_0^T + s \int_0^T e^{-st} f(t)\, dt$$

$$= e^{-sT} f(T) - f(0) + sF(s) \qquad (28)$$

Expression (28) shows a very complex situation for FLT. It shows that a first order differential equation (DE) will lead to a two point boundary value problem. Thus there is some redundancy because it only requires one unknown to solve a first order equation. Note that f(T) in (28) is not known because it depends on the solution. The other two terms in (28) are same as in ILT. This is probably another reason why infinite time was originally considered in the definition of (1), as f(T) term will vanish if T goes to infinity.

Using an example we illustrate how the boundary conditions can be determined to solve a DE. Consider a second order simultaneous DE taken from [1, p. 65]:

$\frac{dy}{dt} = -z, \quad \frac{dz}{dt} = y, \quad y(0) = 1, \quad z(0) = 0$

Taking the FLT on both sides results

$sY(s) - y(0) + e^{-sT}y(T) = -Z(s)$

$sZ(s) - z(0) + e^{-sT}z(T) = Y(s)$

The above can be rewritten in matrix form as

$\begin{bmatrix} s & 1 \\ -1 & s \end{bmatrix} \begin{bmatrix} Y(s) \\ Z(s) \end{bmatrix} = \begin{bmatrix} y(0) - e^{-sT}y(T) \\ z(0) - e^{-sT}z(T) \end{bmatrix}$

Solving them we get

$Y(s) = \frac{s}{s^2+1} - e^{-sT}\frac{s}{s^2+1}y(T) + e^{-sT}\frac{1}{s^2+1}z(T)$ \hfill (29)

$Z(s) = \frac{1}{s^2+1} - \frac{1}{s^2+1}e^{-sT}y(T) - \frac{s}{s^2+1}e^{-sT}z(T)$ \hfill (30)

From the expression for Y(s) in (29) we get

$Y(s) = \frac{1}{s^2+1}[s - e^{-sT}sy(T) + e^{-sT}z(T)]$ \hfill (31)

Since the denominator of (31) is zero at ± i, the numerator of (31) also should be zero at these values to make Y(s) analytic at all points in the s-plane. Therefore we must have, from (31),

$s - e^{-sT}sy(T) + e^{-sT}z(T) = 0, \quad at \ s = \pm i$

Substituting these values of s in the above expression we get the following two equations

$ie^{iT} - iy(T) + z(T) = 0$

$-ie^{-iT} + iy(T) + z(T) = 0$

Solving them simultaneously gives $y(T) = \cos T$, and $z(T) = \sin T$. These values as the boundary conditions in (29) and (30) will help to select the correct solution using (9) and (10).

It should be noted that the solution of a DE is unique. So both ILT and FLT will produce the same t-domain expression over the same finite time interval. However they will produce different expressions in the s-domain. The above application shows that before we try to solve a differential equation we must first solve for the terminal conditions using the analyticity property of the FLT.

## VI  Convolution Theorem

The convolution theorem has been used by the Laplace transform engineers directly and indirectly in many applications. Surprisingly, the FLT literature, [2,4,5] does not talk about the convolution theorem. Therefore the following is a new result in the FLT theory. The convolution theorem is a very fundamental property of the linear time invariant (LTI) systems. The LTI systems are defined by the ordinary differential equations with constant coefficients. This property says that if the impulse response of a LTI system is given by h(t) then the response y(t) for any other input u(t) can be obtained by the following convolution integral:

$$y(t) = \int_0^t u(\tau) h(t-\tau) d\tau \tag{32}$$

It is worth pointing out here that h(t) exists for infinite time and therefore y(t) is also defined for infinite time, independent of the duration for u(t). Using the convolution theorem for ILT [1, p. 92] the above expression (32) can be reduced to

$$Y(s) = H(s).U(s) \tag{33}$$

Thus (33) allows an easy way to find the response of any LTI system. Pictorially (33) implies the block diagrams shown in the Figure 2. The part (b) of the figure, and the convolution theorem, can help us to cascade systems to find the combined output from the input U(s):

$$Y(s) = Y_1(s)H_2(s) = U(s)H_1(s)H_2(s) \tag{34}$$

This configuration is shown in part (c) of Figure 2.

We first show graphically that the convolution theorem cannot be valid for FLT systems. This graphical view point is very simple, basic, and can be found in many engineering text books, for example [16, pp. 63-75]. We normally overlook the point that is important for our subject, so we present the idea in details. Consider the two functions below:

$$f(t) = H_a(T - t), \quad 0 \le t \le T \tag{35}$$

$$g(t) = H_b(T - t), \quad 0 \le t \le T \tag{36}$$

Here H(t) is the Heaviside step function defined by

$$H(t) = \begin{cases} 1, & t \ge 0 \\ 0, & t < 0 \end{cases}$$

These functions are shown in the Figure 3. They are same as (4) only with different notations. The convolution of these two functions is derived as shown below:

$$h(t) = \int_0^t f(\tau) g(t - \tau) d\tau$$

$$= \int_0^t H_a(T - \tau) H_b(T - [t - \tau]) d\tau$$

$$= \int_0^t H_a(T - \tau) H_b(T - t + \tau) d\tau \tag{37}$$

The Figure 3 also shows that the resulting function h(t) from (37) extends beyond the domain of definition of the two functions involved. Because of this reason if we consider FLT for only time T, a portion of h(t) will not be considered in the FLT, and thus the convolution theorem will not work. This will not be apparent from the proof of the theorems for both ILT and the FLT, because both proofs follow almost exactly the same statements as shown below. A closer look will show that the theorem works for the ILT because the ILT assumes infinite time.

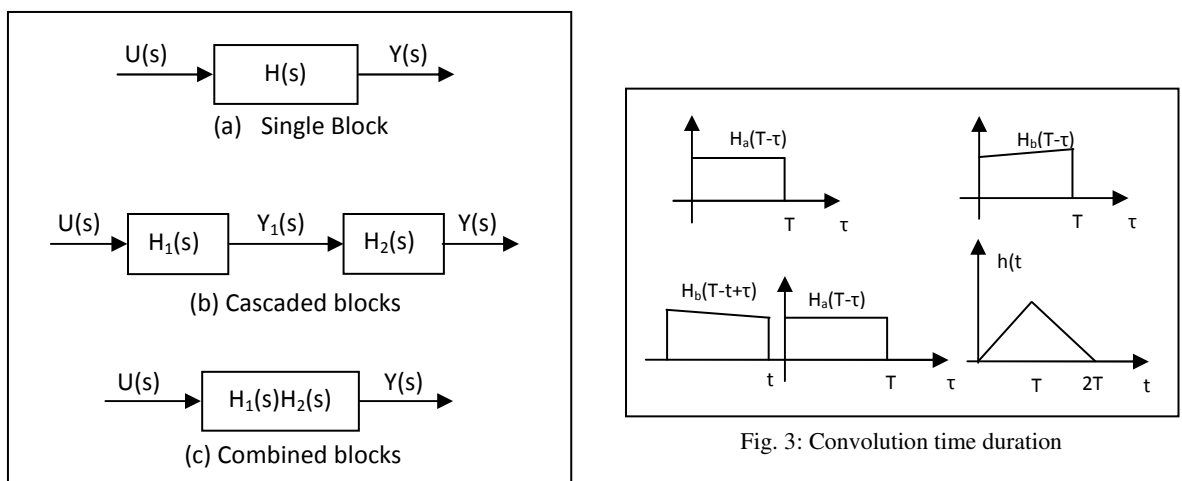

Fig. 2: Use of convolution theorem for ILT

Fig. 3: Convolution time duration

***Theorem:*** If f and g are continuous on $[0, T]$ then the FLT relation given below is true:

$$\mathcal{L}_T[f].\mathcal{L}_T[g] \neq \mathcal{L}_T[f * g] \tag{38}$$

The theorem is proved in the following way [1, pp. 91-93]. Using the definition of FLT (3) we write

$$\mathcal{L}_T[f].\mathcal{L}_T[g] = \left[\int_0^T e^{-sv} f(v)dv\right]\left[\int_0^T e^{-su} g(u)du\right] \tag{39}$$

Since all functions are continuous in the square region bounded by $0 \leq u \leq T$ and $0 \leq v \leq T$ we can move the second integral sign in (39) at the left to get

$$= \int_0^T \left[\int_0^T e^{-s(u+v)} f(v)g(u)du\right] dv \tag{40}$$

Now we define u+v=t and change the limits of the second integral in (40) as shown below

$$= \int_0^T \left[\int_v^{T+v} e^{-st} f(v)g(t-v)dt\right] dv \tag{41}$$

Since the function f(t) is zero for t > T, we can change the upper limit from T+v to T. If t<v the function g(t-v) is zero, so we can set the lower limit v to zero also, without affecting the integral and can write:

$$= \int_0^T \left[\int_0^T e^{-st} f(v)g(t-v)dt\right] dv \tag{42}$$

Again switching the order of integration in (42) we get

$$= \int_0^T \left[\int_0^T e^{-st} f(v)g(t-v)dv\right] dt \tag{43}$$

Now take the exponential outside the inner integral of (43), because the variable of integration is v inside and write:

$$= \int_0^T \left[ e^{-st} \int_0^T f(v)g(t-v)dv \right] dt \qquad (44)$$

Since v cannot be greater than t for the function g(t-v) is zero in that region, therefore the upper limit T for v can be set to t, giving us.

$$= \int_0^T e^{-st} \left[ \int_0^t f(v)g(t-v)dv \right] dt \qquad (45)$$

Now we have the convolution expression inside (45). We have seen from Figure 3 that the convolution integral in (45) goes beyond T and the expression (45) neglects that part of the function between [T, 2T], because the limits on the outer integral in (45) is 0 to T, and therefore (45) cannot be equal to (46) given below:

$$\mathcal{L}_T[f * g] \qquad (46)$$

This concludes the proof of the convolution theorem for the FLT system.

Thus we cannot multiply two FLT transfer functions and use the result for a single constant time interval application, which is normally what we do in ILT engineering. It should be mentioned that if we consider [0, 2T] as the interval for convolution in (45) then the equality will hold. But the time [T, 2T] is not available for the engineering process, because during that time the next activity has to be scheduled by the operating system. Also extending it to 2T will eventually mean to 4T and finally to infinity. The reason that the convolution theorem holds in ILT expressions, is because ILT expressions are by definition valid only for infinite time, as can be seen from (11) and (12). In the following section we show how we can go around the convolution theorem and still use the FLT theory to design engineering systems. Our approach is a direct extension of the conventional ILT method. This approach has not appeared in the literature before, as far as we know.

VII    Control System Design

The Infinite Laplace Transform (ILT) theory is very widely used for the design of various engineering control systems in military and commercial applications. The core idea behind this approach can be described with the help of the block diagram shown in Figure 4. In this figure all variables represent ILT functions. The transfer function P(s) represents the known plant or the process that is to be controlled. C(s) is the controller that needs to be designed. R(s) is the reference input and Y(s) is the output and also provides the negative feedback. The objective of the design is to find the parameters of the controller C(s) so that the error e(t), the inverse ILT of E(s), meets the engineering requirements.

The existing classical design method is basically a trial and error approach. The poles and zeros of the controller C(s) are selected to match the objective. Ideally, of course C(s) could be selected as the inverse of P(s) [17, p. 24] but for many reasons that cannot be done. However, in the absence of convolution theorem, and finite time constraints, that concept is not meaningful any more. In this section we illustrate the basic idea of an FLT based design method using an open loop system and a simple differential equations (DE). The objective is only to demonstrate that an ILT based trial and error method can be extended to the FLT theory. In view of the lack of convolution theorem, we present a slightly different approach. Consider the plant model:

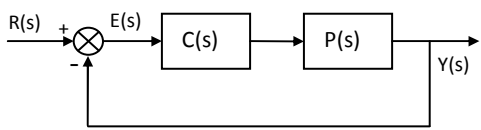

Fig.4: A control system design approach

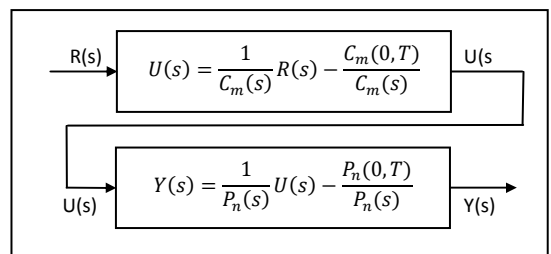

Fig.5: FLT design principle for controller

$$p_n \frac{d^n y}{dt^n} + p_{n-1} \frac{d^{n-1} y}{dt^{n-1}} + \cdots + p_1 \frac{dy}{dt} + p_0 y = u(t) \tag{47}$$

The FLT of n-th order derivative of a function can be written as

$$\mathcal{L}_T[y^{(n)}(t)] = s^n \mathcal{L}_T[y(t)] - \sum_{k=0}^{n-1} s^k\, y^{(n-k-1)}(0) + e^{-sT} \sum_{k=0}^{n-1} s^k\, y^{(n-k-1)}(T)$$

We use the following simplifying notation

$$y_n(0,T) = -\sum_{k=0}^{n-1} s^k\, y^{(n-k-1)}(0) + e^{-sT} \sum_{k=0}^{n-1} s^k\, y^{(n-k-1)}(T)$$

Using the above notation the FLT of the DE for the plant (47) can be written as

$$p_n s^n Y(s) + p_n y_n(0,T) + p_{n-1} s^{n-1} Y(s) + p_{n-1} y_{n-1}(0,T) + \cdots + p_1 s^1 Y(s) + p_1 y_1(0,T) + p_0 Y(s) = U(s)$$

Using more simplifying notations the last expression can be reduced to

$$P_n(s) Y(s) + P_n(0,T) = U(s) \tag{48}$$

In (48) we have borrowed the notation style from the software engineering concepts. The number of parameters in the parentheses will define the meaning of the same variable P. That is

$$P_n(s) = \sum_{k=0}^{n} p_k s^k \quad \text{and} \quad P_n(0,T) = \sum_{k=0}^{n} p_k\, y_k(0,T)$$

Therefore from (48) we can express Y(s) of the plant as:

$$Y(s) = \frac{1}{P_n(s)} U(s) - \frac{P_n(0,T)}{P_n(s)} \tag{49}$$

Similarly the controller model C(s) can be written as in (50):

$$U(s) = \frac{1}{C_m(s)} R(s) - \frac{C_m(0,T)}{C_m(s)} \tag{50}$$

Observe that we do not have the concept of transfer function in (49) or (50). That is because the input U(s) in (49) is only a part of the right hand side, preventing us in taking a ratio of output over input. Even if the initial conditions are zero the right hand side of (49) will still be non-zero, because it includes the final time values. This is also another reason why we cannot multiply two transfer functions in the FLT based systems. We point out again that both expressions (49) and (50) are analytic, i.e. there are no poles. Instead of transfer functions we will use the term system functions for expressions like (49).

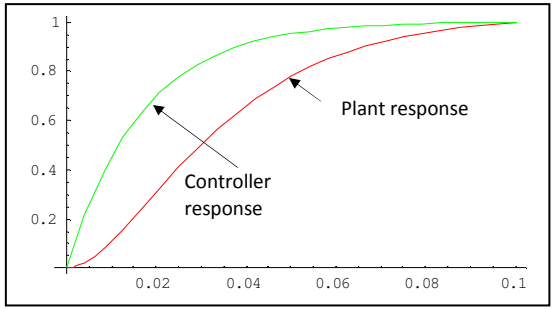

| Controller parameter c=60 | | | | |
|---|---|---|---|---|
| Controller output data | | | | |
| 0.418288 | .662041 | .80409 | .88687 | .935107 |
| .963218 | .979602 | .989147 | .994709 | .997954 |
| Plant output data | | | | |
| .0905053 | .266025 | .44658 | .601706 | .723942 |
| .815463 | .881729 | .928606 | .961183 | .983512 |

Fig. 6: Response graphs for controller c=60

Table 2: FLT controller design data for c=60

The block diagram in Figure 5 shows the composite open loop system. The design objective is to find the unknown coefficients of the controller $C_m(s)$, for the known plant and the known reference input. Since the convolution theorem prevents us to multiply the two blocks together in our design process, we must consider the output of the individual system functions separately and then feed them to the following blocks as shown in Figure 5.

*A Design Example*

More details of the method described above are given using an example and its numerical simulation. Consider a first order known plant, given by the DE

$$\dot{y} + py = u(t), \qquad 0 \leq t \leq T \tag{51}$$

Here u(t) is the output from the controller and is also the input to the plant, and p is a plant parameter of known value. Assume also a first order controller with unknown controller parameter c and a known input r(t) as shown in the following DE

$$\dot{u} + cu = r(t) = 1, \qquad 0 \leq t \leq T \tag{52}$$

For simplicity all initial conditions are assumed to be zero. The FLT of the controller (52) gives

$$sU(s) + e^{-sT}u(T) + cU(s) = \frac{1-e^{-sT}}{s}$$

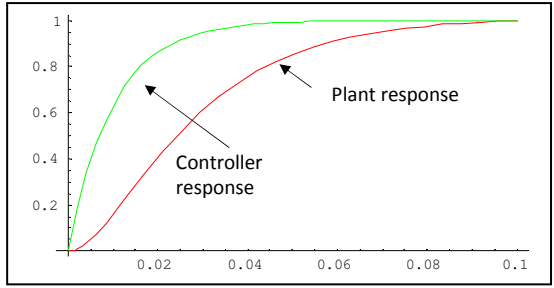

Fig. 7: Response graphs for controller c=100

| Controller parameter c=100 | | | | |
|---|---|---|---|---|
| Controller output data | | | | |
| .593432 | .834586 | .932702 | .972627 | .98877 |
| .99537 | .99811 | .99915 | .999577 | .999825 |
| Plant output data | | | | |
| .133062 | .35704 | .556277 | .706267 | .811346 |
| .882093 | .928693 | .959122 | .978763 | .991314 |

Table 3: FLT controller design data for c=100

$$U(s) = \frac{1}{s+c}\frac{1-e^{-sT}}{s} - \frac{e^{-sT}u(T)}{s+c}$$

$$= \frac{1-e^{-sT}-se^{-sT}u(T)}{s(s+c)} \tag{53}$$

Since U(s) must be analytic we must have, from (53), at s=-c,

$$1 - e^{cT} + ce^{cT}u(T) = 0$$

This gives the boundary value for u(T) as

$$u(T) = \frac{1-e^{-cT}}{c} \tag{54}$$

Now we can similarly write the system function of the plant model (51) to give

$$sY(s) + e^{-sT}y(T) + pY(s) = U(s)$$

Solving for Y(s) we can write

$$Y(s) = \frac{1}{s+p}U(s) - \frac{e^{-sT}y(T)}{s+p}$$

We can now substitute for the input U(s) from (53) to get

$$Y(s) = \frac{1}{s+p}\frac{1-e^{-sT}-se^{-sT}u(T)}{s(s+c)} - \frac{e^{-sT}y(T)}{s+p}$$

$$= \frac{1-e^{-sT}-se^{-sT}u(T)-s(s+c)e^{-sT}y(T)}{s(s+c)(s+p)} \tag{55}$$

Again using analyticity of Y(s) at s=-p, we get from (55):

$$1 - e^{pT} + pe^{pT}u(T) + p(c-p)e^{pT}y(T) = 0$$

This produces the expression for the terminal value of y(T) as

$$y(T) = \frac{1-e^{-pT}-pu(T)}{p(c-p)} \tag{56}$$

We can reconfirm the design problem as - find c of the controller so that the output y(t) from the plant matches the given step input to the controller. In this example we have assumed p=50 for the plant. First we try c=60 and find u(T) and y(T) using (54) and (56). Then we find the input to the plant by inverting the controller FLT (53) using Taylor series method. Finally we invert the plant (55), again using Taylor series, to find the output from the plant. The process can be repeated with another trail value for c. The results are shown in the tables and figures for two controller parameters c=60 and c=100. The graphs for the plant and the controller responses for a step input at the controller are also shown in Figures 6 and 7. The graphs show that a faster response from the controller provides a faster response from the plant. The important feature of the design method is that there was no need to use the concept of poles and zeros of the system. As mentioned many times, the FLT theory does not have poles. The simulation used the Mathematica software tools.

The literature search seems to indicate that the first paper on the application of the FLT to control system was presented by Datko [18]. He used state equations, the optimal control theory, and its solutions, for the linear control problems. Using the analyticity property of the FLT he has shown how the optimal value for the final time can be computed. The control law of course was bang-bang, that is, the control switches between two fixed constant values +1 and -1. He has also used the FLT theory for quadratic optimal control problem to find the terminal conditions. His approach is essentially a discrete approach. He finds a discrete fixed number, the final time or the terminal conditions using the FLT theory. [7] Uses a sequence of FLT to generate ILT to analyze tracking control problem. As mentioned earlier, paper [7] also identifies the location for the proof of the analyticity of the FLT. Rosen [19] has used FLT to solve finite time optimal control problems. He has assumed the control law as a linear combination of exponential functions of FLT variables between two sample intervals. It is a combination of state equation and Laplace approach. In our paper we have shown how the classical ILT based concepts can be extended in a straight forward manner to find continuously varying control functions using the FLT theory and its numerical method.

## VIII  Fourier Transform

The Fourier transform (FT) theory is intimately related to the ILT. Thus a brief note on the finite time aspect of the FT should be mentioned in this FLT paper. We point out that the Fourier series can be time limited but the Fourier transform cannot. In case of Fourier Transforms all the harmonic signals must exist over infinite time interval. The corresponding time signals, like in ILT, must also be of infinite time duration. The FT pair is defined by the following two expressions.

$$F(w) = \int_{-\infty}^{\infty} f(t)e^{-jwt} \, dt \tag{57}$$

$$f(t) = \frac{1}{2\pi}\int_{-\infty}^{\infty} F(w)e^{jwt} \, dw \tag{58}$$

The expression (57) gives the FT and the (58) gives the Inverse FT (IFT).

The expression (57) shows that the time function f(t) must exist for infinite time and the IFT expression (58) shows the same thing also. The IFT also shows that the frequency w must go from -∞ to +∞ over the real numbers and not just over integers.

One way to examine this infinity property of FT is to visualize the example of the delta function. Its FT is 1 for all w. That means all cosine functions that create the delta function have unit amplitude and zero phase. If you draw some of these cosine functions [20, pp. 178-179] you will find that the functions are adding up to create the pulse and becoming zero at all other places. This example shows that all cosine functions must be defined over all time, and the same must be true for the delta function also. That is, the delta function must exist as zero for the entire real line except the place where it is non-zero. This infinite time duration is an important restriction of FT as well as ILT theory for their engineering applications. They do not satisfy the theoretical foundation of engineering practice. Again, this infinite time is not a restriction for the Fourier series though. In Fourier series all sine functions, as well as the original time function, can all be of finite duration [21, pp. 31-32]. However the Fourier series places serious restrictions on the bandwidth as explained below.

The sinusoidal functions of the FT theory are orthogonal only over infinite time. For example, if you consider an interval of one milliseconds time, then the first fundamental will have a frequency of 1 KHz and the second harmonic will have 2 KHz. Thus 1500 Hz is not there in the Fourier series over one millisecond interval and is not orthogonal also. Therefore for finite time problems the Fourier series harmonic functions will go beyond bandwidth very quickly. However if you consider infinite time then all of the above frequencies are there and also are orthogonal. If you consider only finite time then these functions will not be orthogonal to each other and therefore will not form the basis vector set and consequently will not be able to represent the original time function correctly. In addition, if you take only a few of the frequencies then more errors will be created because now you have fewer vectors in your basis set. Thus the FT theory, like ILT, is not the correct tool for finite time engineering systems.

## IX Conclusions

We have examined the effect of changing the infinite time requirement for the classical Laplace transform theory on engineering applications. It was necessary because engineering does not use infinite time. It should be considered as a paradigm shift from infinite time to finite time in engineering. It changes our view of engineering. The most important feature of the Finite Laplace transform (FLT) functions is that they are analytic over the entire complex plane and therefore does not have any poles. This feature eliminates the need for frequency and root locus based approaches. It also does not show any stability problems. It is shown that the FLT does not satisfy convolution theorem and therefore requires a new design approach for control systems. The transfer functions of classical theory become system functions. We have given a simple numerical method for the inverse FLT. More research will be required to generate proper software development tools for engineering design and simulation using the FLT theory.

## Acknowledgements

The author is deeply indebted to Professor Joel L. Schiff of Mathematics department, University of Auckland, New Zealand, for many help he has provided so gracefully and quickly during the early development of this research work. During this initial conversation phase the core idea, the infinite time issue, of this paper was crystallized. Incidentally, he is also the author of reference [1] of this paper.

The author remembers that this problem was first discussed with a very respected friend and colleague, Prof. Kalyan Ray of electrical engineering department at Jadavpur University, Kolkata, India when we were graduate students. It took so many years to find the answer.